\documentclass[11pt]{amsart}

\usepackage{amsmath, amssymb, amsfonts, amsthm, bm}

\newcommand{\epsi}{\varepsilon}
\newcommand{\fhi}{\varphi}
\newcommand{\norm}[1]{\lVert#1\rVert}
\newcommand{\numbersystem}[1]{\mathbb{#1}}
\newcommand{\R}{\numbersystem{R}}
\newcommand{\abs}[1]{\lvert#1\rvert}
\newcommand{\card}[1]{\lvert#1\rvert}

\theoremstyle{plain}
\newtheorem{quotedtheorem}{Theorem}
\newtheorem{theorem}{Theorem}

\begin{document}

\title{A lower bound for the equilateral number of normed spaces}
\author{Konrad J. Swanepoel}
\address{Department of Mathematical Sciences,
    University of South Africa, PO Box 392,
    Pretoria 0003, South Africa}
\thanks{This material is based upon work supported by the South African National Research Foundation under Grant number 2053752.
The second author thanks the DGES grant BFM2003-01297 for financial support.
Parts of this paper was written during a visit of the second author to the Department of Mathematical Sciences, University of South Africa, in January 2006.}
\email{\texttt{swanekj@unisa.ac.za}}
\author{Rafael Villa}
\address{Departamento An\'alisis Matem\'atico, Facultad de Matem\'aticas,
Universidad de Se\-villa, c/Tarfia, S/N, 41012 Sevilla, Spain} \email{\texttt{villa@us.es}}
\subjclass[2000]{Primary 46B04; Secondary 46B20, 52A21, 52C17}
\begin{abstract}
We show that if the Banach-Mazur distance between an $n$\nobreakdash-\hspace{0pt}dimensional normed space $X$ and
$\ell_\infty^n$ is at most $3/2$, then there exist $n+1$ equidistant points in $X$. By a well-known
result of Alon and Milman, this implies that an arbitrary $n$-dimensional normed space admits at
least $e^{c\sqrt{\log n}}$ equidistant points, where $c>0$ is an absolute constant.
We also show that there exist $n$ equidistant points in spaces sufficiently close to $\ell_p^n$, $1<p<\infty$.
\end{abstract}
\maketitle

\section{Notation}
Throughout the paper we use the same symbol $c$ for different absolute positive constants. Let $X$
denote a normed space of finite dimension $\dim X=n$. Let $e(X)$ denote the largest size of an
equilateral set in $X$.
As usual, the space $\ell_p^n$, $1\leq p<\infty$, is defined as $\R^n$ with the norm $\norm{(x_1,x_2,\dots,x_n)}_p=(\sum_{i+1^n}\abs{x_i}^p)^{1/p}$, and $\ell_\infty^n$ is $\R^n$ with the norm $\norm{(x_1,x_2,\dots,x_n)}_\infty=\max_i\abs{x_i}$.
The \emph{Banach-Mazur distance} between two $n$-dimensional normed spaces is defined as $d(X,Y)=\inf\norm{T}\,\norm{T^{-1}}$, where the infimum is taken over all linear, invertible operators $T:X\to Y$.
We say that $X$ is a \emph{$(1+\epsi)$-copy} of $Y$ if $d(X,Y)\leq 1+\epsi$.

\section{The main theorems}
It is conjectured \cite{MR2163782,MR25:1492, MR93h:53012, MR43:1051,
MR97f:52001} that $e(X)\geq n+1$ for all $n$-dimensional normed spaces $X$. This is known for $n\leq 3$ \cite{MR43:1051} but open for $n\geq
4$. It is true for spaces sufficiently close to Euclidean:
\begin{quotedtheorem}[Brass \cite{MR2000i:52012} \& Dekster \cite{MR2001b:52001}]\label{bd}
Let $X$ be an $n$-dimensional normed space with Banach-Mazur distance $d(X,\ell_2^n)\leq
1+\frac{1}{n}$. Then an equilateral set in $X$ of at most $n$ points can be extended to one of
$n+1$ points. In particular, $e(X)\geq n+1$.
\end{quotedtheorem}
Combining this theorem with Gordon's estimate \cite{MR87f:60058} in the Dvoretzky theorem \cite[\S 4]{GM}, the
following general lower bound follows:
\[e(X)\geq c(\log n)^{1/3}.\]
We improve this to the following:
\begin{theorem}\label{main}
For any $n$-dimensional normed space $X$ we have $e(X)\geq e^{c\sqrt{\log n}}$, where $c>0$ is an
absolute constant.
\end{theorem}
The proof works by looking for large subspaces close to either $\ell_2^k$ or $\ell_\infty^k$. To
this end we use the following theorem:
\begin{quotedtheorem}[Alon-Milman \cite{AM}]\label{am}
Let $X$ be an $n$-dimensional normed space. Then for each $\epsi>0$ there exists $c=c(\epsi)>0$
such that $X$ either contains a $(1+\epsi)$-isomorphic copy of $\ell_2^m$ for some $m$ satisfying
$\log\log m\geq\frac{1}{2}\log\log n$ or contains a $(1+\epsi)$-isomorphic copy of $\ell_\infty^k$
for some $k$ satisfying $\log\log k > \frac{1}{2}\log\log n -c$.
\end{quotedtheorem}
If the second case occurs in the above theorem, i.e., if we find a $(1+\epsi)$-isomorphic copy of
$\ell_\infty^k$ in $X$, then we need a result similar to Theorem~\ref{bd} for spaces near
$\ell_\infty^n$. This we provide as follows:
\begin{theorem}\label{th}
Let $X$ be an $n$-dimensional normed space with Banach-Mazur distance $d(X,\ell_\infty^n)\leq 3/2$.
Then $e(X)\geq n+1$.
\end{theorem}
We can then choose $\epsi=1/2$ in Theorem~\ref{am} to obtain $e(X)\geq e^{c\sqrt{\log n}}$
in this case.
On the other hand, if we find a $(1+\frac{1}{2})$-isomorphic copy of $\ell_2^m$ where
$m>e^{\sqrt{\log n}}$, we cannot yet apply Theorem~\ref{bd}. We first have to find a
$k$-dimensional subspace of the $m$-dimensional space that is $(1+\frac{1}{k})$-isomorphic to
$\ell_2^k$. This is provided by the following:
\begin{quotedtheorem}[Milman \cite{M}]
Let $X$ be an $m$-dimensional normed space and $0<\epsi<1$. Then $X$ contains a $k$-dimensional
subspace $Y$ where $k\geq c\epsi^2m/d^2(X,\ell_2^m)$ and $d(Y,\ell_2^k)\leq 1+\epsi$.
\end{quotedtheorem}
See also \cite[Corollary~4.2.2]{GM}. Putting $\epsi=1/k$ and $X$ the $m$-dimensional subspace
that is $3/2$-isomorphic to $\ell_2^m$ into the above theorem, we obtain $k > cm^{1/3}$, and
Theorem~\ref{bd} then gives $e(X)>cm^{1/3}> ce^{(1/3)\sqrt{\log n}}$. To complete the proof of
Theorem~\ref{main}, it only remains to prove Theorem~\ref{th}.

\begin{proof}[Proof of Theorem~\ref{th}]
We use the Brouwer fixed point theorem \cite[\S 14.3]{Browder}, as in Brass' proof of Theorem~\ref{bd}. Without loss of
generality we may assume $X=(\R^n,\norm{\cdot})$ and
\[ \norm{x}\leq\norm{x}_\infty\leq\frac{3}{2}\norm{x}\quad\text{ for all $x\in X$.}\]

Let $I=\{(i,j): 1\le i< j\le n+1\}$, with $\card{I}=n(n+1)/2=N$. For $\epsi=(\epsi_{i,j})_{(i,j)\in I}\in[0,\frac{1}{2}]^N$, let
\begin{align*}
p_1(\epsi)&=(-1,0,\dots 0),\\
p_j(\epsi)&=(\epsi_{1,j},\dots\epsi_{j-1,j},-1,0,\dots 0),\qquad 2\le j\le n-1,\\
p_n(\epsi)&=(\epsi_{1,n},\dots\epsi_{n-1,n},-1),\\
p_{n+1}(\epsi)&=(\epsi_{1,n+1},\dots\epsi_{n,n+1}).
\end{align*}

For $1\le i< j\le n$ we have $\norm{p_i(\epsi)-p_j(\epsi)}_\infty= 1+\epsi_{i,j}$.
Define
$\fhi:[0,1/2]^N\to[0,1/2]^N $ by $\fhi_{i,j}(\epsi) =
1+\epsi_{i,j}-\norm{p_i(\epsi)-p_j(\epsi)},\ 1\le i<j\le n$.
Note that
$$
\fhi_{i,j}(\epsi)\ge 1+\epsi_{i,j}-\norm{p_i-p_j}_\infty=0
$$
and
$$
\fhi_{i,j}(\epsi)\le
1+\epsi_{i,j}-\frac{2}{3}\norm{p_i-p_j}_\infty=\frac{1}{3}(1+\epsi_{i,j})\le\frac{1}{2},
$$
so $\varphi$ is well-defined.
Brouwer now gives the existence of a point $\epsi'=(\epsi'_{i,j})\in[0,1/2]^N$ with $\varphi(\epsi')=\epsi'$, which implies that $\norm{p_i(\epsi')-p_j(\epsi')}=1$ for all $1\le i< j\le m$.
We have obtained $n+1$ equilateral points.
\end{proof}

%

%
%
%
%

\section{A generalization to $\ell_p^n$}
The following theorem partially generalizes Theorem~\ref{th} to all $\ell_p^n$ spaces with
$1<p<\infty$.

\begin{theorem}\label{th2} For each $n>2$ and $p\in(1,\infty)$ there exists $R(p,n)>1$ such that for any $n$-dimensional normed space $X$ with Banach-Mazur
distance $d(X,\ell_p^n)\le R(p,n)$ we have $e(X)\geq n$.
In fact,
\begin{align*}
R(p,n) &= \max_{\theta>0} \left(\frac{1+(1+\theta)^p}{2+(n-2)\theta^p}\right)^{1/p}\\
& \sim 1+\frac{p-1}{2p}n^{-\frac{1}{p-1}}\text{ as $n\to\infty$ with $p$ fixed.}
\end{align*}
\end{theorem}

\begin{proof} We follow the proof of Theorem~\ref{th}.
Assume $X=(\R^n,\norm{\cdot})$ and
$$
\norm{x}\le\norm{x}_p\le R\norm{x}\quad\text{ for all $x\in X$.}
$$
Fix $\beta, \gamma>0$.
Let $I=\{(i,j): 1\le i< j\le n\}$, with $\card{I}={n(n-1)/2}=N$.
For $\epsi=(\epsi_{i,j})_{(i,j)\in I}\in[0,\beta]^N$, let
\begin{align*}
p_1(\epsi)&=(-\gamma,0,\dots 0),\\
p_j(\epsi)&=(\epsi_{1,j},\dots \epsi_{j-1,j},-\gamma,0,\dots 0),\qquad 2\le j\le n-1,\\
p_n(\epsi)&=(\epsi_{1,n},\dots \epsi_{n-1,n},-\gamma).
\end{align*}

For $1\le i< j\le n$ we have
$$
\|p_j-p_i\|_p^p= \sum_{k=1}^{i-1}\abs{\epsi_{k,j}-\epsi_{k,i}}^p+ (\epsi_{i,j}+\gamma)^p+
\sum_{k=i+1}^j \epsi_{k,j}^p+\gamma^p.
$$
Define $\fhi:[0,\beta]^N\to[0,\beta]^N$ by $\fhi_{i,j}(\epsi) =
1+\epsi_{i,j}-\norm{p_i-p_j}$ for $1\le i<j\le n$.
On the one hand,
\begin{align*}
\fhi_{i,j}(\epsi) &\le 1+\epsi_{i,j}-R^{-1}\norm{p_i-p_j}_p
\\
& \le 1+\epsi_{i,j}-R^{-1}\left[(\gamma+\epsi_{i,j})^p+\gamma^p\right]^{1/q}.
\end{align*}
Taking into account that the latter is increasing with respect to $\epsi_{i,j}$, the inequality
$\epsi_{i,j}\le\beta$ implies
$$
\fhi_{i,j}(\epsi) \le 1+\beta-R^{-1}\left[(\gamma+\beta)^p+\gamma^p\right]^{1/q}.
$$
Therefore, if $(\gamma+\beta)^p+\gamma^p\ge R^p$ then $\fhi_{i,j}(\epsi)\le\beta$. On the other hand
\begin{align*}
\fhi_{i,j}(\epsi) &\ge 1+\epsi_{i,j}-\norm{p_i-p_j}_p\\
&\ge 1+\epsi_{i,j}- \left[ (n-2)\beta^p+(\gamma+\epsi_{i,j})^p+\gamma^p \right]^{1/q}.
\end{align*}
Again the latter is increasing with respect to $\epsi_{i,j}$, so using $\epsi_{i,j}\ge0$ we have
$$
\fhi_{i,j}(\epsi) \ge 1- \left[ (n-2)\beta^p+2\gamma^p \right]^{1/q}.
$$
Then $\fhi_{i,j}(\epsi_1,\dots \epsi_m)\ge0$ would follow if $(n-2)\beta^p+2\gamma^p \le 1. $
Subsequently, if
\begin{equation}\label{cond}
(\gamma+\beta)^p+\gamma^p\ge R^p\qquad \text{and}\qquad (n-2)\beta^p+2\gamma^p \le 1,\tag{$\ast$}
\end{equation}
then $\fhi$ is well defined. Brouwer now gives a point
$\epsi'=(\epsi'_{i,j})\in[0,\beta]^N$ such that $\fhi(\epsi')=\epsi'$, implying that the points $p_1(\epsi'),\dots p_n(\epsi')$ are equilateral.

Finally, to take the best choice for the parameters in \eqref{cond} we have to maximize the expression
$(\gamma+\beta)^p+\gamma^p$
under the constraints
$(n-2)\beta^p+2\gamma^p \le 1$ and $\beta , \gamma\ge0$.
Setting $\theta=\beta/ \gamma$, we obtain
$$
R^p=\max_{\theta>0}\frac{1+(1+\theta)^p}{2+(n-2)\theta^p}.
$$
It is not difficult to see that for $\theta$ close to $n^{-1/(p-1)}$ the right-hand side is $>1$ and $R-1\sim \frac{p-1}{2p}n^{-\frac{1}{p-1}}$.
\end{proof}

\section{Concluding remarks}
For $p=2$ the estimate in the above theorem is $d(X,\ell_2^n)\lesssim1+\frac{1}{4n}$, slightly worse than Theorem~\ref{bd}.
However, we don't know how to obtain $n+1$ equidistant points as in Theorem~\ref{th}.
It would also be interesting to know whether arbitrary equilateral sets of at most $n$ points in spaces near $\ell_p^n$ can be extended as in Theorem~\ref{bd}.
A different idea will be needed to extend the above theorem to the case $p=1$.
See \cite{Swanepoel-survey} for a survey on equilateral sets, as well as \cite{MR1995795, Smyth, Swanepoel-AdM} for further results on equilateral sets in $\ell_p^n$.

\end{document}